\documentclass[12pt, reqno]{amsart}
\usepackage{amsmath, amsthm, amscd, amsfonts, amssymb, graphicx, color}
\usepackage{setspace}
\usepackage{mathrsfs}
\usepackage{multicol}
\usepackage[bookmarksnumbered, colorlinks, plainpages]{hyperref}
\hypersetup{colorlinks=true,linkcolor=red, anchorcolor=green, citecolor=cyan, urlcolor=red, filecolor=magenta, pdftoolbar=true}

\newtheorem{theorem}{Theorem}[section]
\newtheorem{lemma}[theorem]{Lemma}

\theoremstyle{definition}

\theoremstyle{remark}

\numberwithin{equation}{section}

\newcommand{\CC}{\mathbb {C}}
\newcommand{\CU}{\mathbb {C}^+}
\newcommand{\Hg}{\mathscr{H}}
\newcommand{\R}{\mathbb{R}}

\newcommand{\D}{\mathbb{D}}
\begin{document}
\setcounter{page}{1}
\title[Volterra type integral and composition  operators on model spaces] {Volterra type integral operators and composition operators on model spaces }
\author [Tesfa  Mengestie]{Tesfa  Mengestie }
\address{Department of Mathematical Sciences \\
Stord/Haugesund University College (HSH)\\
Klingenbergvegen 8, N-5414 Stord, Norway}
\email{Tesfa.Mengestie@hsh.no}
\subjclass[2010]{Primary 47B32,30H20;Secondary 46E22,46E20,47B33 }
\keywords{Volterra type integral  operator,  bounded, compact, composition operator, model spaces}
\begin{abstract}
We  study some mapping properties of Volterra type integral operators and composition operators
on model spaces. We also discuss  and give out a couple of  interesting  open problems  in model spaces where any possible solution of the problems can be  used to  study  a number of other operator theoretic related  problems in the spaces.
\end{abstract}
 \maketitle
\section{Introduction}
We denote by   $\mathscr{H}$ the space of holomorphic functions. We define the Volterra type integral operator on
$\mathscr{H}$ induced by a homomorphic symbol $g$ by
 \begin{equation*}
 V_gf(z)= \int_0^z f(w)g'(w) dw.
\end{equation*}
Questions about various   operator
theoretic properties of $V_g$ expressed in terms of function
theoretic conditions on $g$ have been a subject of high interest
mainly after the works of Pommerenke \cite{Pom} and later by Aleman, Cima and Siskakis \cite{ALC, Alsi1,Alsi2}. In  \cite{Alsi1}, it was proved that  $V_g$    on the Hardy space  $H^p(\D), 1\leq p<\infty,$  is bounded if and only if  $g$ belongs to the space of $BMOA,$ while compactness is characterized  in terms of the corresponding $VMOA$ spaces.

In this paper we plan to study the bounded, compact and Hilbert--Schmidt  properties of  the
Volterra type integral  operators  and composition operators acting between the  model spaces $K_I^2$ and the Hardy space $H^2(\CU)$ of the upper-half plane.  We recall  that the model spaces $K_I^2$,  where $I$ is an inner function, are the subspace of $H^2(\CU)$ defined by
$K_I^2= H^2(\CU)\ominus IH^2(\CU).$ The space $K_I^2$ is a reproducing kernel Hilbert space with kernel function
$$K_w(z)=\frac{i}{2\pi}\cdot\frac{1-\overline{I(w)}I(z)}{z-\overline{w}},\ \ \ \  k_w(z)= K_w(z)/\left\|K_w\right\|_2.$$
Each inner function $I$ has either one or infinitely many components. We call it a one-component if there exists a $\delta$ in $ (0,1)$ for which the level set $\left\{ z\in \CU: |I(z)|< \delta \right\}$ is connected. For model spaces generated by such class of inner functions, we prove the following.
\begin{theorem}\label{thm1}
Let $I$ be a one-component inner function on $\CU$ and $g$ be a holomorphic function on $\CU$. Then $V_g: K_I^2 \to H^2(\CU)$ is
\begin{enumerate}
\item  bounded if and only if
\begin{equation}
\sup_{w\in\CU} \int_{\CU} |k_w(z)|^2 |g'(z)|^2 \Im z  dA(z)<\infty;
\end{equation} \item   compact if and only if
\begin{equation}
\lim_{|w| \to \infty} \int_{\CU} |k_w(z)|^2 |g'(z)|^2 \Im z  dA(z)=0;
\end{equation}
\item  in  the Hilbert--Schmidt class $\mathcal{S}_2(K_I^2, H^2(\CU))$ if and only if
\begin{equation}
\int_{\CU} |g'(z)|^2 (1-|I(z)|^2) dA(z) <\infty,
\end{equation} where $dA$ denotes the usual  Lebesgue area measure.
\end{enumerate}
\end{theorem}
The criteria  found  above  classify the  bounded and compact  $V_g$ on  Model spaces in terms of their action on the reproducing kernels. As will be seen in Section \ref{proof}, by converting the boundary integral in Hardy spaces  to an area integral in terms of derivatives, we could reduce  the  boundedness and compactness problems of $V_g$  into questions of bounded and compact Carleson embedding maps (measures) for $K_I^2.$ Such measures have been identified when $I$ belongs to
 the class of one-component inner functions: for example see  \cite{Al, Al1,Co2, Co1, Co3, TV}. Despite several attempts made by many experts in the subject, a complete identification of the  measures remains  open when $I$ belongs to the class of infinitely many component inner functions.   For special class of model spaces which admit reproducing kernel Riesz bases corresponding to sparse sequence of points, which  includes  the de Branges spaces, the author together  with K. Seip and Y. Belov \cite{BMS2} have described the measures in  terms  of a condition analogues to the $A_2$ weight condition. In what follows we will use those  results in \cite{BMS2} to  study  the bounded and compact properties of $V_g$. We note in passing that those  spaces constitute   model spaces generated by the class of  infinitely many component inner functions.
  \subsection{The Volterra type integral operators on $\Hg(\Gamma, v)$ }
 Let $\Gamma=(\gamma_n)$ be a sequence of
distinct complex numbers and $v=(v_n)$ be a weight sequence that
satisfies the admissibility condition
\begin{equation}
\sum_{n=1}^\infty \frac{v_n}{1+|\gamma_n|^2}<\infty.
\label{adm}
\end{equation}  Any such pair $(\Gamma, v)$ parameterizes the space $\Hg(\Gamma,v)$ which
  consists of all functions
\begin{equation}
\label{new0} f(z)=\sum_{n=1}^\infty\frac{a_n v_n}{z-\gamma_n}
 \end{equation}
for which
\begin{equation}
\label{new1} \|f\|_{\Hg(\Gamma,v)}= \|(a_n)\|_{\ell_v^2}<\infty, \ \  \ell_v^2=
\big\{ (a_n) : \|a\|_v^2=\sum_{n=1}^\infty |a_n|^2v_n
<\infty\big\}\end{equation}
 and $z$ belongs to the set
\begin{equation}
\label{new2}
 (\Gamma,
v)^{*}=\Big\{ z\in \CC:\ \sum_{n=1}^\infty
\frac{v_n}{|z-\gamma_n|^2}<\infty\Big\}.
\end{equation}
The equation in \eqref{new0} means that we obtain the value of a function $f$ in
$\Hg(\Gamma, v)$ at a point $z$ in $(\Gamma, v)^{*}$ by computing
the
 weighted discrete Hilbert transform:
\begin{equation}
 (a_n)\mapsto\sum_{n=1}^\infty \frac{a_n
v_n}{z-\gamma_n}, \label{discreteH}
\end{equation}
which is well defined whenever $(a_n)$ belongs to $\ell_v^2.$ This follows from  an application of Cauchy--Schwarz inequality
along with \eqref{new1} and \eqref{new2}.
   The de Branges
 spaces, model subspaces of the Hardy space  $H^2$ which admit Riesz bases of reproducing kernels
  and the Fock-type spaces
 studied in \cite{BL} are  all examples of spaces of the kind  $\Hg(\Gamma, v).$  We refer to
\cite{BMS2, BMS, TYM} for detailed accounts. The Carleson measures on $\Hg(\Gamma, v)$  have
been described  in \cite{BMS2} when $\Gamma$ grows at least
exponentially, i. e., when
\begin{equation}
\label{expon}
 \inf_{n\geq1}\frac{|\gamma_{n+1}|}{|\gamma_n|}>1.\end{equation}
 When we consider such sparse sequence $\Gamma,$ it is natural
 to partition $\CU$ in the following way. We set
 $\Omega_1= \left\{ z\in \CU: |z| <(|\gamma_1| +|\gamma_2|)/2\right\}$ and then for $n\geq 2$
 $$\Omega_n= \left\{ z\in \CU: (|\gamma_{n-1}|+|\gamma_{n}|)/2 \leq |z| < (|\gamma_{n}|+|\gamma_{n+1}|)/2\right\}.$$
 We now state our first result on $\Hg(\Gamma, v)$.
\begin{theorem} \label{thm2}
Suppose that the sequence $\Gamma$ satisfies the sparseness
condition \eqref{expon} and $g: (\Gamma, v)^{*} \to \CU$ be homomorphic. Then $V_g:\Hg(\Gamma, v)\to H^2(\CU)$ is
\begin{enumerate}
\item bounded  if and only if
\begin{equation}
\label{trivial} \sup_{n\geq1} \int_{\Omega_n} \frac{v_n |g'(z)|^2 \Im z
dA(z)}{|z-\gamma_n|^2}<\infty
\end{equation} and
\begin{equation}
\label{carleson} \sup_{n\geq1} \left( \sum_{l=1}^{n}v_l
\sum_{m=n+1}^\infty \int_{\Omega_m} \frac{|g'(z)|^2 \Im z dA(z)}{|z|^2} +\sum_{l=
n+1}^\infty \frac{v_l}{ |\gamma_l|^2}
\sum_{m=1}^{n}\int_{\Omega_m} |g'(z)|^2 \Im z dA(z)\right)<\infty;
\end{equation}
\item compact if and only if the  "little oh" counterparts of \eqref{trivial} and \eqref{carleson} holds;
\item in  the Hilbert--Schmidt class $\mathcal{S}_2(\Hg(\Gamma, v), H^2(\CU))$ if and only if
\begin{equation}
\label{trivialHSS} \sum_{n=1}^\infty \int_{\Omega_n} \frac{v_n |g'(z)|^2
\Im z }{|z-\gamma_n|^2}dA(z) <\infty
\end{equation}
and
\begin{equation} \label{BHSS}
\sum_{k=1}^\infty v_k
\sum_{k=n+1}^\infty \int_{\Omega_k} \frac{|g'(z)|^2 \Im z dA(z)}{|z|^2} +\sum_{k=
n+1}^\infty \frac{v_n}{ |\gamma_n|^2}
\sum_{k=1}^{n}\int_{\Omega_k} |g'(z)|^2 \Im z dA(z)<\infty.
\end{equation}
\end{enumerate}
\end{theorem}
Condition \eqref{trivial} of the theorem is a condition about the local
behavior of the symbol $g$, while condition \eqref{carleson} deals with its global
behavior. Combining the two conditions, we see that \eqref{trivial} may be replaced by the stronger global condition
\begin{equation*}
\sup_{n\geq1} \int_{\CU} \frac{v_n |g'(z)|^2 \Im z
dA(z)}{|z-\gamma_n|^2}<\infty.
\end{equation*} In some special cases, this condition, or the original
one \eqref{trivial} automatically gives condition \eqref{carleson}. For example if
the weight sequence $v_n$ grows at least exponentially and the numbers
$v_n/|\gamma_n|^2$ decay at least exponentially with $n$. The same  conclusion follows if
the  weight sequence $v_n$ is also  summable.
\section{The composition operator on $\Hg(\Gamma, v)$  }
Let $\psi$ be a homomorphic function in a given domain. We denote by $C_\psi$ the composition operator $f\longmapsto f\circ\psi.$  As a consequence of the Littlewood's Subordination Theorem \cite{LW}, it has long been known that all composition  operators are bounded on all the Hardy spaces $H^p$ of the unit disc where $0<p\leq \infty$. Although corresponding Hardy spaces of the disc
and half-plane are isomorphic, composition operators act differently in the two  domains. Unlike the case of $H^p(\D)$, not all composition operators are bounded on  $H^p(\CU)$.

  The study of compact composition operators on $H^2(\D)$ first  appeared in H. Schwartz \cite{SCH} thesis in the late sixties though a complete function theoretic characterization,  in terms of the inducing map's Nevanlinnan counting function,  was obtained by Shapiro later in \cite{Shapiro}. Unlike $H^p(\D)$ again, it  was proved in   \cite{VM} that there exists no compact composition operator on the Hardy space $H^p(\CU)$.   The  work in \cite{Shapiro} was continued by Y. Lyubarskii and E. Malinnikova \cite{LE}   for  $C_\psi: K_I^2(\D) \to H^2(\mathbb{D})$  and a complete extension was made when the generating inner function $I$ is one-component. A more general trace ideal criteria for $C_\psi$ was obtained latter  in \cite{ALM}.   The   problem remains open when $I$ has  infinitely many components.

 In this section, we study the bounded, compact and Hilbert--Schmidt  properties of composition operators  on model spaces when the spaces admit  normalized reproducing kernel Riesz bases associated to  sparse sequences of points  in $\CU$. We note that such spaces represent some classes of model spaces
 generated by infinitely many component inner functions.
\begin{theorem} \label{thm3}
Suppose that the sequence $\Gamma$ satisfies the sparseness
condition \eqref{expon} and $\psi$ be a nonconstant analytic function on $\CU$, and $m$ refers to the usual Lebesgue measure on the real line. Then  $C_\psi:\Hg(\Gamma, v)\to H^2(\CU)$ is
\begin{enumerate}
\item bounded  if and only if
\begin{equation}
\label{triviall} \sup_{n\geq1} \sup_{y\geq0} \int_{\psi(\Omega_n)} \frac{v_n
dx}{|\psi(x+iy)-\gamma_n|^2}<\infty
\end{equation} and
\begin{equation}
\label{carlesonn} \sup_{n\geq1}  \left( \sum_{l=1}^{n}v_l
\sum_{k=n+1}^\infty \sup_{y\geq0} \int_{\psi(\Omega_k)} \frac{ dx}{|\psi(x+iy)|^2} +\sum_{l=
n+1}^\infty \frac{v_l}{ |\gamma_l|^2}
\sup_{y\geq0}\sum_{k=1}^{n}  m(\psi(\Omega_k))\right)<\infty,
\end{equation}
where $\psi(\Omega_k)= \{\psi(x+iy): x+iy\in \Omega_k \}$;
\item compact if and only if the "little oh" counterparts of conditions \eqref{triviall} and \eqref{carlesonn} hold;
\item  in the Hilbert--Schmidt class $\mathcal{S}_2(\Hg(\Gamma, v), H^2(\CU))$ if and only if
\begin{equation}
\label{trivialHS} \sum_{k=1}^\infty \sup_{y>0}\int_{\psi(\Omega_k)} \frac{v_k
dx}{|\psi(x+iy)-\gamma_k|^2} <\infty
\end{equation}
and
\begin{equation} \label{BHS}
\sum_{k=1}^\infty v_k \sum_{j=k+1}^\infty \sup_{y>0}\int_{\psi(\Omega_j)}
\frac{dx}{|\psi(x+iy)|^2}+  \sum_{k=1}^\infty
\frac{v_k}{|\gamma_k|^2} \sup_{y> 0}\sum_{j=1}^{k-1} m(\psi(\Omega_j))<\infty.
\end{equation}
\end{enumerate}
\end{theorem}
\section{Open problems}
\subsection{Open problem 1}
As will be seen in the next section, the proof of our boundedness and compactness results for  both $C_\psi$ and $V_g$  rely mainly on  previously obtained  Carleson measure results on $\Hg(\Gamma, v)$.  When  the sequence of points
$\gamma_n= n$ and the weight sequences $v_n \simeq 1$, the space $\Hg(\Gamma, v)$
becomes the classical Paley--Wiener space for which its Carleson measures have long been identified. When the sequences $\gamma_n$ grows at least exponentially, the measures are completely described in terms of $A_2$ type condition in  \cite{BMS2}. The case when $\gamma_n$   grows  between these  two extreme cases remains an interesting  open  problem. A solution of this will settle the long-standing open problem of identifying the Carleson measures for  Model spaces, of course modulo to the existence of reproducing kernel Riesz basis in Model spaces. The novelty and the core of
the approach in \cite{BMS2}  has been in turning a number of problems including  the Carleson measure problem  into questions about
different mapping properties of Hilbert transforms in weighted spaces of functions and
sequences. To extend the approach, one may need to look into all possible interplays between the smootheness (regularity) of the weight
 sequence $(v_n)$ and  the growth of the sequence $\Gamma= (\gamma_n)$.
\subsection{Open problem 2}  Both of our results in parts iii) of Theorem 1.2 and Theorem 2.1 deal with  the Hilbert-Schmidt properties of
the operators  $V_g$ and $C_\psi$ on $\Hg(\Gamma, v)$. A more natural and interesting question is to ask the Schatten $\mathcal{S}_p$ class
membership of these maps  from $\Hg(\Gamma, v)$ into $H^2(\CU)$ for  all $p$ in the range $ 0< p<\infty.$  Dealing with  these problems would have been easier had we been known a complete characterization of the $\mathcal{S}_p$ class  membership of the imbedding maps induced by Carleson measures for the spaces, which itself is still another open problem.
 A good starting point for all of these problems could be to consider first the case when  the sequences
$\gamma_n$ grows at least exponentially. A very special case of the  $\mathcal{S}_p$ membership problem  pertaining to the embedding maps may be seen in \cite{TYM1}.

A  word on notation: the  notation $U(z)\lesssim V(z)$ (or
equivalently $V(z)\gtrsim U(z)$) means that there is a constant
$C$ such that $U(z)\leq CV(z)$ holds for all $z$ in the set of a
question. We write $U(z)\simeq V(z)$ if both $U(z)\lesssim V(z)$
and $V(z)\lesssim U(z)$.
\section{Proof of the results} \label{proof}
\textit{Proof of Theorem~\ref{thm1}}.
The key tool in the proofs of the first two results is the  Littlewood--Paley
description of Hardy spaces in terms of derivative :
\begin{equation}
\label{LP}
\left\|f\right\|_2^2= 2 \int_{\CU} |f'(z)|^2 \Im z dA(z)
\end{equation} for each $f$ in $ H^2(\CU)$. The formula helps to convert  a boundary integral to an area integral in terms  of derivatives.
    Since $(V_gf)' (z)= f(z) g'(z),$  an application of the above  identity ensures that for  each $f$ in $ K_I^2$ we have
    \begin{equation}
    \label{new3}
    \left\|V_gf\right\|_2^2=2\int_{\CU} |f(z)|^2 |g'(z)|^2 \Im z dA(z).
\end{equation}  It follows from this that the boundedness of  the operator
$V_g: K_I^2 \to H^2(\CU)$ is equivalent to saying that the measure  $\mu_g$  where
\begin{equation}
\label{measure}d\mu_g(z)= |g'(z)|^2 \Im z  dA(z)
\end{equation}  is a Carleson measure for $K_I^2.$
Assuming that $I$ is a one-component inner function, a result of Cohn \cite{Co1}
ensures that the measure property holds if and only if the normalized reproducing kernels of $K_I^2$ are  uniformly embedded into $L^2(\CU, d\mu_g) $ from which the boundedness part of the  result follows.

To prove (ii), first from a standard normal family argument, we have that $V_g$ is compact if and only if  $V_g f_n \to 0$ in $H^2(\CU)$ for each sequence of functions $f_n$ in $ K_I^2$  such that the sequence is uniformly bounded in norm and  converges uniformly to zero on compact subsets of $\CU.$ It means that the compactness of $V_g$ is equivalent to saying that the embedding map from $K_I^2$ into $L^2(\mu_g)$ is compact.  Since $I$ is one-component, then the desired conclusion follows from a result of  Cohn \cite{Co1} again.\\

To prove part (iii), recall that the operator $V_g$ belongs to the Hilbert--Schmidt class if and only if for any orthonormal basis $(e_n)$ of $K_I^2$, the sequence  $(\| V_g e_n \|_{2}^2)$  is summable.  It follows by \eqref{new3} that
 \begin{eqnarray*}
 \sum_{n=1}^\infty\| V_g e_n \|_{2}^2  &\simeq& \sum_{n=1}^\infty \int_{\CU} |e_n(z)|^2 |g'(z)|^2 \Im z dA(z) \\
 &=& \int_{\CU} \sum_{n=1}^\infty|e_n(z)|^2 |g'(z)|^2 \Im z dA(z).
 \end{eqnarray*}  On the other hand, because of the reproducing property of the kernel and Parseval's  identity, we have
 \begin{align}
\label{kernel}
K_{w}(z)=\sum_{n=1}^\infty \langle K_w, e_n\rangle e_n(z) = \sum_{n=1}^\infty e_n(z) \overline{e_n(w)} \  \text{and}\ \ \| K_w\|_2^2= \sum_{n=1}^\infty |e_n(w)|^2.
\end{align} The desired conclusion easily  follows since
   \begin{equation}
  \sum_{n=1}^\infty|e_n(z)|^2 = \|K_z\|_2^2 \simeq \frac{1-|I(z)|^2}{\Im z}.
  \end{equation}
 \textit{Proof of Theorem~\ref{thm2}}.  Applying the Littlewood--Paley identity again, the estimate
\begin{equation}
  \left\|V_g f\right\|_2^2= 2 \int_{\CU} |f(z)|^2 |g'(z)|^2 \Im z dA(z) \lesssim \left\|f\right\|_{\Hg(\Gamma,v)}^2
 \end{equation} holds  if and only if   the embedding map  from $\Hg(\Gamma,v)$ into $ L^2(\CC, \mu_g)$ is bounded where
 $\mu_g$ is the nonnegative measure defined by \eqref{measure}.
 Then, the  boundedness  part of the result immediately follows from an application of
 Theorem~1.1 in \cite{BMS2}.

  By a result of Nordgren's \cite{Nord}, a sequence of functions $f_n$ in $\Hg(\Gamma,v)$ converges weekly to zero if and only if it converges pointwise to zero and  $\sup_n \left\|f_n\right\|_{\Hg(\Gamma,v)} <\infty.$ It  follows that $V_\psi$ is compact if and only if  the  measure $\mu_g$ is a vanishing Carleson measure for $\Hg(\Gamma,v).$  Then an application of   Theorem~1.2 in  \cite{TYM} finishes the proof of the compactness part.

 It  remains to prove the statement about Hilbert--Schmidt membership. First, we may observe that the reproducing kernel of $\Hg(\Gamma,v)$ at a point
$\lambda$ in $(\Gamma, v)^*$  is explicitly given by
 \[k_\lambda(z)=
\sum_{n=1}^\infty
\frac{v_n}{(\overline{\lambda}-\overline{\gamma_n})(z-\gamma_n)};\]
this is a direct consequence of the definition of the space $\Hg(\Gamma,v)$. Furthermore, the sequence $e_n(z)= \sqrt{v_n}/(z-\gamma_n)$ constitutes
an orthonromal basis to $\Hg(\Gamma,v)$. Then
\begin{eqnarray}
\label{HS}
\|V_g\|_{\mathcal{S}_2}^2& = &\sum_{n=1}^\infty \| V_g e_n\|_{2}^2 \simeq \sum_{n=1}^\infty \int_{\CU} |e_n(z)|^2 |g'(z)|^2 \Im z dA(z)\nonumber\\
&=& \int_{\CU} \|K_z\|_{\Hg(\Gamma,v)}^2  |g'(z)|^2 \Im z dA(z)\nonumber\\
&=& \int_{\CU} \sum_{n=1}^\infty
\frac{v_n}{|z-\gamma_n|^2} |g'(z)|^2 \Im z dA(z).
\end{eqnarray}
Using the growth condition \eqref{expon}, we observe that if $z\in \Omega_m$, then
\begin{equation}
\sum_{n=1}^\infty
\frac{v_n}{|z-\gamma_n|^2} \simeq  \sum_{n=1}^{m-1}
\frac{v_n}{|z|^2}+ \frac{v_m}{|z-\gamma_m|^2} + \sum_{n=m+1}^\infty
\frac{v_n}{|\gamma_n|^2}.
\end{equation} Plugging this estimate in \eqref{HS}, we find that $\|V_g\|_{\mathcal{S}_2}$ is finite if and only if
both \eqref{trivialHSS} and \eqref{BHSS} hold.

\textit{Proof of Theorem~\ref{thm3}}.
Suppose now that $\psi$ is analytic on $\CU$. Then for each fixed $y>0,$ we  define
$\psi_y(x)= \psi(x+iy)$ from $\R$ into $\CU$. We denote by  $m\psi_y^{-1} $ the  pullback measure
\begin{equation*}
m\psi_y^{-1}(E)= m(\psi_y^{-1}(E))
\end{equation*} for each Borel subset $E\subset \CU$ where  $m$ refers to the usual Lebesgue measure on the real line.\\
 Applying change of variables, we observe
\begin{equation}
\int_{\R} |C_\psi f(x+iy)|^2 dx=   \int_{\CU} |f(z)|^2 dm\psi_y^{-1}(z);
\end{equation} which holds because of a classical   result in measure theory (\cite{HL}, p.163 ). This implies that the estimate
\begin{equation}
\label{change}
\|C_\psi f\|_2^2 =  \sup_{y> 0}\int_{\CU}|(f\circ\psi)(x+iy)|^2dx \lesssim \left\|f\right\|_{\Hg(\Gamma,v)}^2
\end{equation}
 holds if and only if  $(m\psi_y^{-1})_{y>0}$ constitute a family of  Carleson measure for $\Hg(\Gamma, v)$ with a uniform bound on $y$ in the sense that
$$\sup_{y>0} \frac{\int_{\CU} |f(z)|^2 dm\psi_y^{-1}(z) }{\left\|f\right\|_{\Hg(\Gamma,v)}^2}<\infty.$$
Then an application of  Theorem~1.1 of  \cite{BMS2} gives  the first part of the result
in the theorem.

The compactness part of our result follows  easily by similar arguments used to prove the
compactness part in  Theorem~\ref{thm2} along with Theorem~1.2 of \cite{TYM}. Thus,  we remain to prove part (iii) of the result. Proceeding
as in the proof of part (iii) of Theorem~\ref{thm2} and applying \eqref{change} we have
\[
\|C_\psi \|_{\mathcal{S}_2}^2= \sum_{n=1}^\infty \| C_\psi e_n\|_{2}^2 =\sum_{n=1}^\infty  \sup_{y> 0}\int_{\CU}|(e_n\circ\psi)(x+iy)|^2dx\]
\begin{equation*}=\sum_{n=1}^\infty  \sup_{y> 0}\int_{\CU} |e_n(z)|^2 dm\psi_y^{-1}(z)= \sum_{n=1}^\infty  \sup_{y> 0}\int_{\CU}\frac{v_n}{|z-\gamma_n|^2}dm\psi_y^{-1}(z)\end{equation*}
\begin{equation}
=  \sum_{n=1}^\infty  \sup_{y> 0}\sum_{j=1}^\infty \int_{\Omega_j}\frac{v_n}{|z-\gamma_n|^2}dm\psi_y^{-1}(z)
\end{equation}
Applying the growth condition in \eqref{expon}, we observe that if $z\in \Omega_n$, then
\begin{align}
\label{a}
\sum_{n=1}^\infty v_n \sup_{y> 0}\sum_{j=1}^\infty \int_{\Omega_j}\frac{dm\psi_y^{-1}(z)}{|z-\gamma_n|^2}\simeq
\sum_{n=1}^\infty v_n \sup_{y> 0}\sum_{j=1}^{n-1} \int_{\Omega_j}\frac{dm\psi_y^{-1}(z)}{|z|^2}\\
+\sum_{n=1}^\infty   \sup_{y> 0}\Bigg( \int_{\Omega_n}\frac{ v_n dm\psi_y^{-1}(z)}{|z-\gamma_n|^2}
+\sum_{j=n+1}^{\infty} v_n\int_{\Omega_j}
\frac{dm\psi_y^{-1}(z)}{|\gamma_n|^2}\Bigg).
\label{aa}
\end{align}
The sum on the right hand side of  \eqref{a} can be further rewritten as
\begin{equation}
\sum_{n=1}^\infty v_n \sup_{y> 0}\sum_{j=1}^{n-1} \int_{\Omega_j}\frac{dm\psi_y^{-1}(z)}{|z|^2}=
\sum_{n=1}^\infty v_n \sup_{y> 0}\sum_{j=1}^{n-1} \int_{\psi(\Omega_j)}\frac{dx}{|\psi(x+iy)|^2}.
\label{b}
\end{equation}
On the other hand,  the sum on the right hand side of \eqref{aa} becomes
\begin{align}
\sum_{n=1}^\infty v_n  \sup_{y> 0}\Bigg( \int_{\Omega_n}\frac{dm\psi_y^{-1}(z)}{|z-\gamma_n|^2}
+\sum_{j=n+1}^{\infty} \int_{\Omega_j}
\frac{v_n}{|\gamma_n|^2}dm\psi_y^{-1}(z)\Bigg)\nonumber\\
 \simeq \sum_{n=1}^\infty  v_n\sup_{y> 0} \int_{\psi(\Omega_n)}\frac{dx}{|\psi(x+iy)-\gamma_n|^2}\nonumber\\
+ \sum_{n=1}^\infty \frac{v_n}{|\gamma_n|^2}  \sup_{y> 0} \sum_{j= n+1}^\infty  m (\psi(\Omega_j)),
\label{bb}
\end{align}
where $m (\psi(\Omega_j))= m\Big(\{\psi(x+iy): x+iy\in \Omega_j\}\Big)$.
From the above series of estimates, we see that $\|C_\psi \|_{\mathcal{S}_2}$ is finite if and only if  the sums in \eqref{b}
and \eqref{bb} are finite, and from which our assertion follows.

%The approach used in \cite{BMS2} may still allow a moderate weaking of the growth of $\gamma_n$ provided that the weight sequences $v_n$ has  certain regularity property. In general it remains a challenge.

\end{document}